\documentclass[11pt]{article}
\usepackage{amsmath,amssymb,a4wide}
\usepackage{color}
 \def\Aut{\mathop{\rm Aut}\nolimits}

 \def\Gal{\mathop{\rm Gal}\nolimits}

 \def\Spec{\mathop{\rm Spec}\nolimits}
 
 \def\deg{\mathop{\rm deg}\nolimits}

\def\Quot{\mathop{\rm Quot}\nolimits}

\def\ord{\mathop{\rm ord}\nolimits}

\def\GL{\mathop{\rm GL}\nolimits}

\def\sep{{\rm sep}}

\let\phi\varphi
\let\epsilon\varepsilon

\def\sp{{\rm spec}}

\newtheorem{Thm}{Theorem}[section]
\newtheorem{Prop}[Thm]{Proposition}

\newtheorem{Cor}[Thm]{Corollary}

\newtheorem{Def}[Thm]{Definition}

\newtheorem{Rem}[Thm]{Remark}

\newtheorem{Ques}[Thm]{Question}

\def\qed{{\hskip0pt\unskip\unskip\nobreak\hfil\penalty50
          \hskip1em\hbox{}\nobreak\hfil
           {$\square$}
          \parfillskip=0pt\finalhyphendemerits=0
          \par}\medskip}

\newenvironment{Proof}
               {\noindent{\bf Proof.}\ }
               {\qed}

               {\noindent{\bf Proof of #1.}\ }
               {\qed}

\newcommand{\BF}{{\mathbb{F}}}

\newcommand{\BQ}{{\mathbb{Q}}}

\newcommand{\Fp}{{\mathfrak{p}}}

\newcommand{\FP}{{\mathfrak{P}}}

\newcommand{\CO}{{\cal O}}

\newbox\mybox
\def\arrover#1{\mathrel{
       \setbox\mybox=\hbox spread 1.4em
              {\hfil$\scriptstyle#1$\hfil}
       \vbox{\offinterlineskip\copy\mybox
             \hbox to\wd\mybox{\rightarrowfill}}}}

\def\larrover#1{\mathrel{
       \setbox\mybox=\hbox spread 1.4em
              {\hfil$\scriptstyle#1\vphantom{g}$\hfil}
       \vbox{\offinterlineskip\copy\mybox
             \hbox to\wd\mybox{\leftarrowfill}}}}

\def\ontoover#1{\mathrel{
       \setbox\mybox=\hbox spread 1.4em
              {\hfil$\scriptstyle#1\vphantom{g}$\hfil}
       \vbox{\offinterlineskip\copy\mybox
             \hbox to\wd\mybox{\rightarrowfill\hskip-2.8mm
                               $\rightarrow$}}}}
\def\leftontoover#1{\mathrel{
       \setbox\mybox=\hbox spread 1.4em
              {\hfil$\scriptstyle#1\vphantom{g}$\hfil}
       \vbox{\offinterlineskip\copy\mybox
             \hbox to\wd\mybox{$\leftarrow$\hskip-2.8mm
                               \leftarrowfill}}}}
\let\longto\longrightarrow

\begin{document}

\title{Drinfeld modular polynomials in higher rank II : \\ Kronecker congruences}
\author{Florian Breuer\footnote{Supported by the Alexander von Humboldt foundation, and by the NRF grant BS2008100900027.}
\and Hans-Georg R\"uck}
\maketitle

\begin{abstract}
This is a sequel to the paper \cite{BR09}, in which we introduced Drinfeld modular polynomials of higher rank, using an analytic construction. These polynomials relate the isomorphism invariants of Drinfeld $\BF_q[T]$-modules of rank $r\geq 2$ linked by isogenies of a specified type.  In the current paper, we give an algebraic construction of greater generality, and prove a generalization of the Kronecker congruences relations, which describe what happens when modular polynomials associated to $P$-isogenies are reduced modulo a prime $P$. We also correct an error in \cite{BR09}.
\end{abstract}


\section{Isomorphism invariants}


Denote by $\BF_q$ the finite field of $q$ elements.
Let $A=\BF_q[T]$ and $r\geq 2$ a positive integer. Let $g_1,\ldots,g_{r-1}$ be algebraically independent over $k=\BF_q(T)$, and let $B=A[g_1,\ldots,g_{r-1}] = \BF_{q}[T,g_1,\ldots,g_{r-1}]$. Denote by $\tau : x\mapsto x^q$ the $q$-Frobenius, and for any ring $R$ we denote by $R\{\tau\}$ the ring of twisted polynomials in $\tau$ with coefficients in $R$, subject to the commutation relations $\tau a = a^q\tau$ for all $a\in R$.

Let $\phi : A \to B\{\tau\}$ be the Drinfeld $A$-module of rank $r$ of generic characteristic  determined by
\[
\phi_T = T + g_1\tau + \cdots + g_{r-1}\tau^{r-1} + \tau^r.
\]
We think of $\phi$ as the {\em monic generic Drinfeld $A$-module of rank $r$.}
A general reference for Drinfeld modules is \cite[\S4]{GossBS}.

The group $\BF_{q^r}^*$ acts on $B\otimes_{\BF_q}\BF_{q^r}$ via 
\[
\lambda * (g_i \otimes 1) = g_i \otimes \lambda^{q^i-1},\qquad i=1,2,\ldots,r-1,
\]
and we denote by 
\[
C := (B\otimes_{\BF_q}\BF_{q^r})^{\BF_{q^r}^*}\cap B 
\]
the subring of invariants with coefficients in $\BF_q$. Equivalently, $C$ consists of those polynomials in $B$ whose monomials are of the form $f(T)g_1^{e_1}g_2^{e_2}\cdots g_{r-1}^{e_{r-1}}$ satisfying $f(T)\in A=\BF_q[T]$ and 
\[
\sum_{k=1}^{r-1} e_k(q^k-1) \equiv 0 \bmod q^r-1. 
\]

Recall that an {\em $A$-field} is a field $L$ together with a ring homomorphism $\gamma : A\to L$.

\begin{Prop}\label{invariants}
For every algebraically closed $A$-field $L$, there is a canonical bijection between the set of isomorphism classes of rank $r$ Drinfeld $A$-modules over $L$, and ring homomorphisms $m : C\to L$ satisfying $m|_A = \gamma$.
\end{Prop}

\begin{Proof}

Let $\psi:A \to L\{\tau\}$ be a rank $r$ Drinfeld $A$-module over the algebraically closed $A$-field $L$, then up to isomorphism we may assume that it is monic, i.e. that
\[
\psi_T = \gamma(T) + a_1\tau  + \cdots + a_{r-1}\tau^{r-1} + \tau^r,\qquad a_i\in L.
\]
We associate to $\psi$ the ring homomorphism
\[
m_\psi:B \rightarrow L \quad \text{ with $m_\psi(g_i) = a_i$  and $m|_A=\gamma$.}
\]
Let $\iota:C \hookrightarrow B$ be the inclusion, then the homomorphism
\[
m_\psi \circ \iota:C \rightarrow L
\]
is invariant under isomorphisms of the Drinfeld module $\psi$.

Conversely, let $m : C \to L$ be a ring homomorphism  with $m|_A=\gamma$. 
For each $i=1,2,\ldots,r-1$ we have $g_i^{q^r-1}\in C$, therefore 
we can extend $m$ to $B$ by defining $m(g_i)$ to be a chosen $(q^r-1)$st root of $m(g_i^{q^r-1})$ in $L$.

Each homomorphism $m:B \rightarrow L$ yields a Drinfeld module $\psi$ with
\[
\psi_T =  \gamma(T) + m(g_1)\tau + \cdots + m(g_{r-1})\tau^{r-1} + \tau^r.
\]

Now let $\psi$ and $\tilde{\psi}$ be two Drinfeld modules such that
$m_\psi \circ \iota = m_{\tilde{\psi}} \circ \iota$.
We have to show that $\psi$ and $\tilde{\psi}$ are isomorphic over $L$.

The group $\BF_{q^r}^*$ acts on $B\otimes_{\BF_q}L = (B\otimes_{\BF_q}\BF_{q^r})\otimes_{\BF_{q^r}} L$ by
$\lambda * (g_i\otimes l) = g_i\otimes \lambda^{q^i-1}l$ and $m_{\psi}$ extends to $m_{\psi} : B\otimes_{\BF_q}L \to L$ as
\[
m_{\psi}(g_i\otimes l) = m_{\psi}(g_i)l.
\]

Suppose that $\psi$ and $\tilde{\psi}$ are not isomorphic, then for each
$\lambda \in \BF_{q^r}^*$, there is a $g_\lambda \in
\{g_1,g_2\ldots,g_{r-1}\}$ for which $m_\psi(\lambda * (g_\lambda\otimes 1)) \neq
m_{\tilde{\psi}}(g_\lambda\otimes 1)$. 
We consider the element $f \in B\otimes_{\BF_q} L$ defined by
\[
f = 1 \otimes 1 - \prod_{\lambda \in \BF_{q^r}^*} \frac{g_\lambda \otimes 1 -
1 \otimes m_\psi(\lambda * (g_\lambda\otimes 1))}{1 \otimes \big(m_{\tilde{\psi}}(g_\lambda\otimes 1)
- m_\psi(\lambda * (g_\lambda\otimes 1))\big)}\,.
\] 
We evaluate
\[
m_{\tilde{\psi}}(f) = 0 \text{ and } m_\psi(\mu * f) = 1 \text{ for each } \mu \in \BF_{q^r}^* \,.
\]
This yields
\[
1 = m_\psi(\prod_{\mu \in \BF_{q^r}^*}(\mu * f)) =  (m_\psi \circ \iota)(\prod_{\mu \in \BF_{q^r}^*}
(\mu * f)) = (m_{\tilde{\psi}} \circ \iota)(\prod_{\mu \in \BF_{q^r}^*}(\mu * f)) = 0 \,,
\]
a contradiction, which proves the proposition.
\end{Proof}

Let $\psi$ be a rank $r$ Drinfeld $A$-module over the $A$-field $L$ given by
\[
\psi_T = \gamma(T) + a_1\tau + \cdots + a_r\tau^r, \quad \text{$a_1,a_2,\ldots,a_r\in L$ and $a_r\neq 0.$}
\]
Denote by $\bar{L}$ the algebraic closure of $L$ and let $\delta\in \bar{L}$ be a $(q^r-1)$st root of $a_r$. Then $\psi$ is isomorphic (over $\bar{L}$) to $\psi' := \delta\psi\delta^{-1}$, given by
\[
\psi'_T = \gamma(T) + a'_1\tau + \cdots + a'_{r-1}\tau^{r-1} + \tau^r, \quad \text{$a'_k = \delta^{1-q^k}a_k$ for $k=1,2,\ldots r-1$.}
\] 

Let $\displaystyle J = \sum_{i=1}^n f_i(T)g_1^{e_{i,1}}g_2^{e_{i,2}}\cdots g_{r-1}^{e_{i,r-1}} \in C$ be an invariant and set
\begin{eqnarray*}
J(\psi) & := & \sum_{i=1}^n \gamma\big(f_i(T)\big)(a'_1)^{e_{i,1}}(a'_2)^{e_{i,2}}\cdots (a'_{r-1})^{e_{i,r-1}} \\
& = & \sum_{i=1}^n \gamma\big(f_i(T)\big)a_1^{e_{i,1}}a_2^{e_{i,2}}\cdots a_{r-1}^{e_{i,r-1}} \delta^{\sum_{i=1}^n\sum_{k=1}^{r-1}e_{i,k}(1-q^k)} \\
& = & \sum_{i=1}^n \gamma\big(f_i(T)\big)a_1^{e_{i,1}}a_2^{e_{i,2}}\cdots a_{r-1}^{e_{i,r-1}} a_r^{-\sum_{i=1}^n\frac{1}{q^r-1}\sum_{k=1}^{r-1}e_{i,k}(q^k-1)} \; \in L.
\end{eqnarray*}

\begin{Cor}\label{isom}
Let $\psi$ and $\tilde{\psi}$ be two rank $r$ Drinfeld $A$-modules over the $A$-field $L$. Then $\psi$ and $\tilde{\psi}$ are isomorphic over the algebraic closure $\bar{L}$ if and only if $J(\psi)=J(\tilde\psi)$ for all $J\in C$.
\end{Cor}

\begin{Proof}
Let $m_\psi : C\to \bar{L}$ be the homomorphism associated to the isomorphism class of $\psi$ by Proposition \ref{invariants}. Then $m_\psi(J) = J(\psi)$ for all $J\in C$ and the result follows.
\end{Proof}

\begin{Rem}
In fancier language, Proposition \ref{invariants} states that $\Spec(C)$ over $\Spec(A)$ is the coarse moduli scheme parametrizing isomorphism classes of rank $r$ Drinfeld $A$-modules over algebraically closed $A$-fields. This was first shown by I.Y. Potemine in \cite{Potemine}.
\end{Rem}

%
%
%
%

If $r=2$, then in fact $C=A[j]$, where  $j=g_1^{g+1} = g^{q+1}/\Delta$ is the usual $j$-invariant of  $\psi_T = \gamma(T) + g\tau + \Delta\tau^2$. 

In general, the ring $C$ is a finitely generated $A$-algebra, and an explicit set of generators is constructed in \cite{Potemine}. This means that one only  needs a finite set of invariants in order to determine whether or not any two Drinfeld modules are isomorphic. For a given finite set of Drinfeld modules of generic characteristic, however, one can find a single invariant to distinguish between them.

\begin{Prop}\label{Distinguish}
Let $S$ be a finite set of pairwise non-isomorphic rank $r$ Drinfeld $A$-modules  of generic characteristic. Then there exists $J\in C$ such that $J$ {\em distinguishes} $S$, i.e. $J(\phi_1)\neq J(\phi_2)$ for all $\phi_1\neq \phi_2\in S$. 
\end{Prop}

\begin{Proof}
We use induction on $n:=|S|$. The result is clearly true for $n\leq 2$. Let
$n\geq 3$, and suppose that the statement is true for  sets of cardinality $n-1$. Pick $\phi_1\in
S$ and let $S_1:=S\smallsetminus\{\phi_1\}$. By the induction hypothesis,
there exists $J_1\in C$ which distinguishes $S_1$. If $J_1$ distinguishes $S$
then we're done. If not, then there exists $\phi_2\in S$ such that
$J_1(\phi_1)= J_1(\phi_2)$, and moreover $\{\phi_1,\phi_2\}$ is
the only pair in $S$ on which $J_1$ takes the same value.

Pick $J_2\in C$ which distinguishes $\{\phi_1,\phi_2\}$, and consider $J_a:=J_1+aJ_2$ for all $a\in A$. If $J_a$ distinguishes $S$ for some $a\in A$, then we're done. If not, then there exists a pair $\psi_1,\psi_2\in S$ for which $J_a(\psi_1)=J_a(\psi_2)$ for at least two distinct values of $a\in A$. From this, we easily deduce that $J_1(\psi_1)=J_1(\psi_2)$, which forces $\{\psi_1,\psi_2\}=\{\phi_1,\phi_2\}$, and $J_2(\psi_1)=J_2(\psi_2)$, which is a contradiction.
\end{Proof}

\begin{Rem}
The conclusion of Proposition \ref{Distinguish} can fail in special characteristic. For example, let $L=\BF_2$ of characteristic $\ker\gamma = T\BF_2[T]$ and $r=3$. Then the three Drinfeld $\BF_2[T]$-modules defined over $L$ by 
\[
\psi^1_T = \tau^3,\qquad \psi^2_T = \tau + \tau^3, \qquad \psi^3_T = \tau^2+\tau^3
\]
are pair-wise non-isomorphic, as witnessed by the invariants $g_1^7,\; g_2^7\in C$, but no single $J\in C$ can distinguish between all three since $J(\psi^i)\in\BF_2$ for every $J\in C$ and $i=1,2,3$. 
\end{Rem}

\section{Isogenies and modular polynomials}


Let $N\in A$ be monic, then $\phi_N$ applied to a variable $X$ defines an
$\BF_q$-linear polynomial $\phi_N(X)\in B[X]$ which is monic and separable
over $B$ of degree $q^{r\deg N}$. Let $K = \Quot(B) =
\BF_{q}(T,g_1,\ldots,g_{r-1})$, denote by $K_N$ the splitting field of
$\phi_N(X)$ over $K$, and let $R_N$ be the integral closure of $B$ in $K_N$.
Then the set $\phi[N]\subset R_N$ of roots of $\phi_N(X)$ is an $A$-module
isomorphic to $(A/NA)^r$.
%
%

The Galois group $\Gal(K_N/K)$ respects this $A$-module structure, and in fact it is shown in \cite{B11} that 
\[
\Gal(K_N/K) \cong \Aut(\phi[N]) \cong \GL_r(A/NA).
\]

Let $f$ be an isogeny from $\phi$ to another Drinfeld module $\phi^{(f)}$,
defined over an extension of $K$. This means that
\begin{equation}\label{eq:isog}
f\phi_T = \phi_T^{(f)} f.
\end{equation}
Write $f = f_0 + f_1\tau + \cdots + f_d\tau^d.$
Replacing $f$ by $f_d^{-1}f$ gives
\[
(f_d^{-1}f)\phi_T = (f_d^{-1}\phi_T^{(f)}f_d) \, (f_d^{-1}f),
\]
so if we replace $\phi^{(f)}$ by an isomorphic Drinfeld module, we may assume that $f_d=1$, i.e. $f$ is monic. In this case, there exists a monic $N\in A$ such that
\[
\ker f \subset \phi[N],
\]
and so $f\in R_N\{\tau\}$.
Now comparing the coefficients of the highest powers of $\tau$ in (\ref{eq:isog}) shows that $\phi_T^{(f)}\in R_N\{\tau\}$ is also monic. For any invariant $J\in C$ we find that $J(\phi^{(f)})\in R_N.$

\begin{Def} Any isogeny $f$ of $\phi$ satisfying $f\subset \phi[N]$ is called an {\em $N$-isogeny}.
\end{Def}

Denote by $I_N$ the set of all monic $N$-isogenies $f\in R_N\{\tau\}$,
as above. Since any such $f$ is
determined by its kernel, the set $I_N$ is finite.

\begin{Def}
Let $J\in C$ be an invariant. We call
\[
\Phi_{J,N}(X) := \prod_{f\in I_N}\big(X - J(\phi^{(f)})\big) \in R_N[X]
\]
the {\em full modular polynomial of level $N$ associated to $J$}.
\end{Def}

%
%

\begin{Def}
Let $H\subset\phi[N]\cong (A/NA)^r$ be an $A$-submodule. Any isogeny $f\in
I_N$ for which $\ker f$ is an element of the $\GL_r(A/NA)$ orbit
of $H$ is called an isogeny of {\em type} $H$.

Let $J\in C$ be an invariant. Then we call
\[
\Phi_{J,H}(X):=\prod_{\text{$f\in I_N$ of type $H$}}\big(X-J(\phi^{(f)})\big)  \in R_N[X]
\]
the {\em modular polynomial of type $H$ associated to $J$.}
\end{Def}

\begin{Prop}
Let $H\subset\phi[N]$ be an $A$-submodule,
and $J\in C$. Then $\Phi_{J,H}(X)\in C[X]$. 
Furthermore, if $J\in C$ distinguishes the Drinfeld modules $\{\phi^{(f)} \;|\; \text{$f$ of type $H$}\}$, then $\Phi_{J,H}(X)$ is irreducible over~$K$.
\end{Prop}

\begin{Proof}
By construction, $\Phi_{J,H}(X)\in R_N[X]$ 
and $\Gal(K_N/K)$ permutes the roots of $\Phi_{J,H}(X)$, so its coefficients lie in $R_N\cap K = B$.

We next show that $\Phi_{J,H}(X)\in C[X]$.
Let $\lambda\in\BF_{q^r}^*$ and consider the isomorphic Drinfeld module $\psi := \lambda^{-1}\phi\lambda$, so 
\[
\psi_T = T + a_1\tau + \cdots + a_{r-1}\tau^{r-1} + \tau^r,
\] 
with each $a_i = \lambda^{q^i-1}g_i \in \BF_{q^r}[T,g_1,\ldots,g_{r-1}] \cong B\otimes_{\BF_q}\BF_{q^r}.$

We note that $\psi[N] = \lambda^{-1}(\phi[N])$, and so every $N$-isogeny $f : \phi \to \phi^{(f)}$ of type $H\subset \phi[N]$ corresponds to an $N$-isogeny $\lambda^{-1} f \lambda : \psi \to \psi^{(\lambda^{-1} f \lambda)} = \lambda^{-1} \phi^{(f)} \lambda$ of type $\lambda^{-1} H\subset \psi[N]$.

Thus, if we specialise the coefficients of $\Phi_{J,H}(X)$ via $g_i\mapsto a_i = \lambda^{q^i-1}g_i$, we obtain the monic polynomial $\Phi_{J(\psi),H}(X)$ whose roots are precisely the $J(\psi^{(\lambda^{-1} f\lambda)})$ for $f\in I_N$ of type $H$. But $J(\psi^{(\lambda^{-1} f\lambda)}) = J(\lambda^{-1} \phi^{(f)} \lambda) = J(\phi^{(f)})$, so $\Phi_{J(\psi),H}(X) = \Phi_{J,H}(X)$ and we have shown that the coefficients of $\Phi_{J,H}(X)$ are $\BF_{q^r}^*$-invariant in $B\otimes_{\BF_q}\BF_{q^r}$. Thus $\Phi_{J,H}(X)\in C[X]$.

Lastly, the condition on $J$ ensures a one-to-one correspondence between the roots of $\Phi_{J,H}(X)$ and $\{\phi^{(f)} \;|\; \text{$f$ of type $H$}\}$.
Since $\Gal(K_N/K)\cong\GL_r(A/NA)$ acts transitively on these roots, $\Phi_{J,H}(X)$ is irreducible over $K$.
\end{Proof}

\begin{Cor}
The full modular polynomial satisfies $\Phi_{J,N}(X)\in C[X]$, and if $J$ distinguishes the Drinfeld modules $\{\phi^{(f)} \;|\; f\in I_N\}$, then the $K$-irreducible factors of $\Phi_{J,N}(X)$ are precisely the $\Phi_{J,H}(X)$ for various types $H\subset \phi[N]$. \qed
\end{Cor}

When $r=2$, then $C=A[j]$ and $\Phi_{j,(A/NA)}(X) = \Phi_N(j,X)\in A[j,X]$ is the usual Drinfeld modular polynomial constructed in \cite{Bae}.

Let $\psi$ be a rank $r$ Drinfeld $A$-module defined over an $A$-field $L$. Then the coefficients of $\Phi_{J,H}(X)\in C[X]$ can be applied to $\psi$ (equivalently, $m_\psi$ from Proposition \ref{invariants} can be applied to each coefficient), resulting in a polynomial 
\[
\Phi_{J(\psi),H}(X) \in L[X], 
\]
whose roots are precisely the $J(\psi^{(f)})$ for $N$-isogenies $f : \psi \to \psi^{(f)}$ of type $H$, which we will prove below. If $N$ is not prime to the characteristic $\ker(\gamma)$ of $\psi$, then we need to make the notion of ``$N$-isogeny of type $H$'' more precise.


For our purposes we define a level-$N$ structure on $\psi$ to be a surjective $A$-module homomorphism 
\[
\mu : \phi[N] \longto \psi[N]\subset \bar{L}. 
\]
When $N$ is prime to the characteristic $\ker(\gamma)$ of $\psi$, then $\mu$ is an isomorphism. If $\mu$ and $\mu'$ are two level-$N$ structures on $\psi$, then there is a $\sigma\in\GL_r(A/NA)$ such that $\mu' = \mu\circ\sigma$.

Let $H\subset \phi[N]$ be an $A$-submodule. Then an isogeny $f : \psi \to \psi^{(f)}$ is said to be of type $H$ if $f(X) = \prod_{h\in H'}(X-\mu(h))$, where $H'\subset \phi[N]$ is an element of the $\GL_r(A/NA)$-orbit of $H$. The set of isogenies of $\psi$ of type $H$ is independent of the chosen level structure $\mu$. We have

\begin{Prop}\label{eval}
Let $\psi$ be a rank $r$ Drinfeld $A$-module over the $A$-field $L$, let $J\in C$ be an invariant and $H\subset \phi[N]$ an $A$-submodule. Then
\[
\Phi_{J(\psi),H}(X) = \prod_{\text{$f$ of type $H$}}\big(X-J(\psi^{(f)})\big) \in L[X].
\]
\end{Prop}

\begin{Proof}
We have
\[
\phi_T(X) = \prod_{u\in\phi[T]}(X - u) = TX + g_1X^q + \cdots + g_{r-1}X^{q^{r-1}} + X^{q^r},
\]
so the $g_i \in \BF_q\big[u \; : \; u\in \phi[T]\big]$ are polynomials over $\BF_q$ in the $u$'s. Moreover, these polynomials are invariant under the $\GL_r(A/TA)$-action on $\phi[T]$.

Next, let $f : \phi \to \phi^{(f)}$ be an isogeny of type $H$, then for a suitable $H'\subset \phi[N]$ we have
\[
f(X) = \prod_{h\in H'}(X-h) = f_0X + f_1X^q + \cdots + f_{d-1}X^{q^{d-1}} + X^{q^d},
\]
where $d = \dim_{\BF_q}(H)$. Again, we see that the $f_i \in \BF_q[h \; : \; h\in H']$ are polynomials over $\BF_q$ in the $h$'s.

Write $\phi^{(f)} = T + g'_1\tau + \cdots + g'_{r-1}\tau^{r-1} + \tau^r$, then comparing coefficients of $\tau^{dr-1}, \tau^{dr-2}, \ldots, \tau^{dr-r+1}$ in 
$f\cdot \phi_T = \phi^{(f)}_T \cdot f$, we obtain 
\[
g'_i \in \BF_q[T,g_1,\ldots,g_{r-1},f_1,\ldots,f_{d-1}] \subset \BF_q[u,h \; : \; u\in\phi[T],\; h\in H'].
\]

As a result 
\[
J(\phi^{(f)}) \in \BF_q[u,h \; : \; u\in\phi[T],\; h\in H'],
\]
but $J \in C\subset \BF_q[T,g_1,\ldots,g_{r-1}]$, so in fact $J(\phi^{(f)}) \in \BF_q\big[u \; : \; u\in\phi[T]\big]$.

Now, replacing $\psi$ by an isomorphic Drinfeld module over $\bar{L}$ if necessary, we may assume that 
\[
\psi_T = \gamma(T) + a_1\tau + \cdots + a_{r-1}\tau^{r-1} + \tau^r, \quad a_1,\ldots,a_{r-1}\in\bar{L},
\]
is monic. Let $M=\mathrm{lcm}(T,N)$ and let $\mu : \phi[M] \to \psi[M]$ be a level-$M$ structure on $\psi$. By the same arguments as above, replacing each $u$ by $\mu(u)$ and each $h$ by $\mu(h)$, we find that each
\[
J(\psi^{(f)}) \in \BF_q\big[\mu(u) \; : \; u\in\phi[T]\big]
\]
is the same polynomial as $J(\phi^{(f)})$, but with each $u$ replaced by $\mu(u)$.

Applying the map $\mu$ to the polynomials $g_i\in \BF_q\big[u \; : \; u\in\phi[T]\big]$, each $g_i$ is mapped to $a_i$, hence $\mu$ coincides there with the homomorphism $m_\psi : C \to \bar{L}$ from Proposition \ref{invariants}. Thus, for each isogeny $f$ of type $H$, we obtain
\[
m_\psi\big(J(\phi^{(f)})\big) = J(\psi^{(f)}),
\]
which completes the proof.
\end{Proof}

\section{Correction to \cite{BR09}}

In \cite{BR09} we gave an analytic construction of modular polynomials of type $(A/NA)^{r-1}$. These polynomials also classify {\em incoming} isogenies $\phi'\to\phi$ with kernels isomorphic to $A/NA$, whereas the point of view of the present article is to classify modular polynomials by the kernels of the dual, {\em outgoing} isogenies $f : \phi\to\phi'$, which explains the shift in terminology from type $A/NA$ to type $(A/NA)^{r-1}$.

Theorem 1.1 of \cite{BR09} claims that the polynomials $\Phi_{J,(A/NA)^{r-1}}(X)$ are irreducible, but this is only true if $J\in C$ distinguishes the set of Drinfeld modules $\{\phi^{(f)} \;|\; \text{$f$ of type $(A/NA)^{r-1}$}\}$, i.e. when $\Phi_{J,(A/NA)^{r-1}}(X)$ has distinct roots. Such invariants $J\in C$ always exist, by Proposition~\ref{Distinguish}.

\section{Reduction mod $P$ and Kronecker congruence relations}

In this section, we let $N=P\in A$ be a monic prime, and we study the reduction of modular polynomials modulo $P$. When we reduce polynomials in $R_P[X]$, then we are actually reducing modulo a chosen prime of $R_P$ extending $PB$ (remember that $R_P$ is integral over $B$), but we will still write mod $P$ for ease of notation.

Define $\BF_P := A/PA$.
We start with the following 
basic result.


\begin{Prop}\label{ordinary}  We have
\[\phi_P \equiv \tilde{\phi}_P \cdot\tau^{\deg(P)} \bmod P \,,\] where
$\tilde{\phi}_P\in (B\otimes_A\BF_P)\{\tau\}$ is not divisible by $\tau$. 

In other words, $\phi$ has ordinary reduction at every prime $P$.
\end{Prop}

We will give two proofs of this result, starting with a conceptual proof.

\medskip

\begin{Proof}
Denote by $\bar{\phi} : A \to (B\otimes_A \BF_P)\{\tau\}$ the reduction of $\phi$ modulo $P$. We have $\bar\phi[P]\cong (A/PA)^{r-h}$, where $h$ is the height of $\bar\phi$ (see \cite[\S4.5]{GossBS}).
The linear term of $\phi_P$ is $P$, so $\bar{\phi}_P$ is divisible by $\tau$ and thus $h\geq 1$. 
We will show that $h=1$ by constructing a specialization of $\phi$ with ordinary reduction at a prime above $P$.

Recall that $k=\BF_q(T)$ and let $F/k$ be a separable extension of degree $r$ which has only one place above the place of $k$ with uniformizer $1/T$ (i.e. $F/k$ is {\em purely imaginary}) and in which $P$ splits completely (such a field exists, by \cite[Chapter X, Theorem 6]{ArtinTate}). Denote by $R$ the integral closure of $A$ in $F$, and let $\psi$ be a rank 1 Drinfeld $R$-module, which is automatically a rank $r$ Drinfeld $A$-module with complex multiplication by $R$ (for example, let $\psi$ be the Drinfeld module corresponding to the lattice $R$ in the algebraic closure of $\BF_q((\frac{1}{T}))\big)$. Let $PR = \Fp_1\Fp_2\cdots \Fp_r$ be the factorization of $PR$ in $R$. Let $L/F$ be a finite extension over which $\psi$ is defined, let $\FP_1$ be a place of $L$ above $\Fp_1$, and denote by $\CO_{\FP_1}$ the valuation ring of $\FP_1$. Since all rank 1 Drinfeld modules have potential good reduction by \cite[Cor. 4.10.4]{GossBS}, we may write
\[
\psi_T = T + a_1\tau + a_2\tau^2 + \cdots + a_r\tau^r, \quad\text{with}\quad a_1,\ldots,a_{r-1}\in\CO_{\FP_1}, \; a_r\in\CO_{\FP_1}^*.
\]
After possibly replacing $L$ by a finite extension and $\psi$ by an isomorphic Drinfeld module, we may assume furthermore that $a_r=1$, so that $\psi$ is the image of $\phi$ under the specialization $g_i\mapsto a_i, \; i=1,2,\ldots, r-1.$

Denote by $\bar\psi$ the reduction of $\psi$ modulo $\FP_1$, then the $P$-torsion module of $\bar\psi$ is 
\begin{eqnarray*}
\bar\psi[P] & = & \bar\psi[PR] \;\cong\; \bar\psi[\Fp_1]\times\bar\psi[\Fp_2]\times\cdots\times\bar\psi[\Fp_r] \\
& \cong & \{0\} \times (A/PA)^{r-1},
\end{eqnarray*}
as required.
%
\end{Proof}

We also provide an elementary proof of Proposition \ref{ordinary}.

\medskip

\begin{Proof}
We specialize $\phi$ to $\psi$, where $\psi_T = T + g_1\tau + \tau^r$, and show that $\psi$ has ordinary reduction at $P$.

For each positive integer $i$, we study the $3^i$ monomials $m=m(T,g_1,\tau)$ in $(\psi_T)^i$, which are constructed from the non-commuting factors $T, g_1\tau$ and $\tau^r$. We call $m$ {\em simple} if no $\tau^r$-factor was used in the construction of $m$. The monomial $m$ can be written uniquely in the form $T^a g_1^b \tau^c$ and we set $\deg_{g_1}(m) := b$ and $\deg_{\tau}(m) := c$.

By induction on $i$, one readily shows that 
\[
\deg_{g_1}(m) \leq \frac{q^{\deg_\tau(m)}-1}{q-1},
\]
with equality if and only if $m$ is simple.

Now suppose $P=\sum_{i=0}^s a_iT^i,$ with $a_s=1$, is our monic prime in $A$ of degree $s$. We consider the coefficient $b_s(T,g_1)$ of $\tau^s$ in
\[
\psi_P = \sum_{i=0}^s a_i(\psi_T)^i = \sum_{j=0}^{rs} b_j(T,g_1) \tau^j.
\]
One of the terms of $b_s(T,g_1)$ is $g_1^{(q^s-1)/(q-1)}$, arising from the monomial $(g_1\tau)^s$ in $(\psi_T)^s$, whereas all other terms have strictly smaller degree in $g_1$, since they arise from non-simple monomials in $(\psi_T)^i$ with $i\leq s$. It follows that $b_s(T,g_1)$ does not vanish modulo $P$, and so the reduction of $\psi$ modulo $P$ has height 1, as required.
\end{Proof}

\begin{Cor}\label{Frobenius}
There exists a unique monic $P$-isogeny $f_0$ of $\phi$ satisfying
\begin{enumerate}
	\item $f_0 \equiv \tau^{\deg P} \bmod P$,
	\item $U_0 := \ker f \cong A/PA = \BF_P$ as $A$-modules, and
	\item $f_0  \phi = \phi^{(f_0)}  f_0$, where
	\[
	\phi^{(f_0)}_T \equiv T + g_1^{|P|}\tau + \cdots + g_{r-1}^{|P|}\tau^{r-1} + \tau^r \bmod P.
	\]
\end{enumerate}
\end{Cor}

\begin{Proof}
Let $U_0 = \ker\big(\phi[P]\longto \bar\phi[P]\big)$ denote the kernel of reduction modulo $P$. Then by Proposition \ref{ordinary}, $U_0\cong A/PA$, and
\[
f_0(X) := \prod_{u_0\in U_0}(X-u_0) \equiv X^{|P|} \bmod P.
\]
It is now easy to verify that $f_0$ has all the required properties.
\end{Proof}

Since the kernel of any $P$-isogeny $f$ is an $\BF_P$-vector space, its kernel $U:=\ker f$ satisfies either $U_0\cap U = \{0\}$, in which case we call $f$ {\em ordinary}, or else $U_0\subset U$, in which case we call $f$ {\em special}. Equivalently, $f$ is ordinary if the reduction of $f$ modulo $P$ is separable, and special otherwise.

\begin{Def}
Let $H\subset \phi[P]$ be an $A$-submodule and $J\in C$ an invariant. We define the following factors of the modular polynomial $\Phi_{J,H}(X)$:

\begin{eqnarray*}
\Phi^{\sp}_{J,H}(X) & := & \prod_{\text{$f\in I_P$ special of type $H$}} \big(X - J(\phi^{(f)})\big) \in R_P[X],\quad \text{and} \\
\Phi^{\ord}_{J,H}(X) & := & \prod_{\text{$f\in I_P$ ordinary of type $H$}} \big(X - J(\phi^{(f)})\big) \in R_P[X].
\end{eqnarray*}
\end{Def}
Clearly $\Phi_{J,H}(X) = \Phi^{\sp}_{J,H}(X)\cdot \Phi^{\ord}_{J,H}(X)$.

We are now ready to prove our main result.

\begin{Thm}[Kronecker Congruence Relations]\label{Kronecker}
Let $P\in A$ be a monic prime, $J\in C$ an invariant and $1\leq s < r$. Then
\begin{enumerate}
\item $\displaystyle \Phi^{\ord}_{J,(A/PA)^s}(X) \equiv \left(\Phi^{\sp}_{J,(A/PA)^{s+1}}(X^{|P|})\right)^{|P|^{s-1}} \bmod P$, and
\item $\displaystyle \Phi_{J,(A/PA)^s}(X) \equiv  \Phi^\sp_{J,(A/PA)^{s}}(X) \cdot \left( \Phi_{J, (A/PA)^{s+1}}^\sp(X^{|P|})\right)^{|P|^{s-1}} \bmod P$.
\end{enumerate}
Furthermore,
\begin{enumerate}
\setcounter{enumi}{2}
\item The reductions of $\Phi_{J,(A/PA)^s}^\sp(X)$ and $\Phi_{J,(A/PA)^{s}}^{\ord}(X)$ modulo $P$ lie in \mbox{$(C\otimes_A\BF_P)[X]$} for every $s=1,\ldots,r$.
\end{enumerate}
\end{Thm}

\begin{Proof}
Let $1\leq s < r$.
Let $f_U$ and $f_{\tilde{U}}$ be two ordinary $P$-isogenies of type $H\cong (A/PA)^s$ with kernels $U$ and $\tilde{U}$, respectively.  
We call $f_U$ and $f_{\tilde{U}}$ {\em equivalent} if $U+U_0 = \tilde{U}+U_0$. This way each equivalence class contains $|P|^s$ elements, since the kernel of each isogeny in the equivalence class of $f_U$ is obtained by adding elements of $U_0$ to each of the $s$ basis vectors of $U$. The $P$-isogeny $f_{U+U_0}$ with kernel $U+U_0$ is then special of type $(A/PA)^{s+1}$, and moreover each special isogeny of type $(A/PA)^{s+1}$ arises from an equivalence class of ordinary isogenies of type $(A/PA)^s$ in this way.

We see that
\[
f_{U+U_0}(X) = \prod_{u\in U} \prod_{u_0\in U_0} \big(X - u - u_0\big) \equiv \prod_{u\in U} (X-u)^{|P|} \equiv f_U(X)^{|P|} \equiv \tau^{\deg P} \big(f_U(X)\big)\bmod P.
\]
Thus, for the corresponding isogenous Drinfeld modules,
\[
\phi^{(U+U_0)}\cdot\tau^{\deg P}\cdot f_U \equiv \phi^{(U+U_0)}\cdot f_{U+U_0} = f_{U+U_0}\cdot\phi \equiv \tau^{\deg P} \cdot f_U \cdot \phi = \tau^{\deg P}\cdot \phi^{(U)} \cdot f_U \bmod P, 
\]
and we obtain
\[
\phi^{(U+U_0)}\cdot\tau^{\deg(P)} \equiv \tau^{\deg(P)}\cdot\phi^{(U)} \bmod P.
\]
For any invariant $J\in C$ we now have
\[
J(\phi^{(U+U_0)}) \equiv J(\phi^{(U)})^{|P|} \bmod P.
\]

If we combine these results, we calculate
\begin{eqnarray*}
\left(\Phi^{\ord}_{J,(A/PA)^{s}}(X)\right)^{|P|} & \equiv & \prod_{\text{$f_U$ ordinary of type $(A/PA)^{s}$}}\big(X^{|P|} - J(\phi^{(U)})^{|P|}\big) \bmod P \\
& \equiv & \prod_{\text{$f_{(U+U_0)}$ special of type $(A/PA)^{s+1}$}}\big(X^{|P|} - J(\phi^{(U+U_0)})\big)^{|P|^{s}} \bmod P \\
& \equiv & \left(\Phi^\sp_{J,(A/PA)^{s+1}}(X^{|P|})\right)^{|P|^{s}} \bmod P,
\end{eqnarray*}
from which (1.) follows.

Next, (2.) follows from (1.) and $\Phi_{J,(A/PA)^s}(X) = \Phi^\sp_{J,(A/PA)^s}(X)\cdot \Phi^{\ord}_{J,(A/PA)^s}(X)$.

It remains to prove (3.). If $s=r$, then 
\[
\Phi_{J,(A/PA)^r}(X) = \Phi^{\sp}_{J,(A/PA)^r}(X) = X-J,
\] 
since the only isogeny of type $(A/PA)^r$ is the endomorphism $\phi_N$, and of course 
\[
\Phi^{\ord}_{J,(A/PA)^r}(X) = 1.
\]
Now suppose that $\Phi^{\sp}_{J,(A/PA)^s}(X) \bmod P \in (C\otimes_A \BF_P)[X]$ for some $1\leq s\leq r$. Since also $\Phi_{J,(A/PA)^{s-1}}(X) \bmod P \in (C\otimes_A \BF_P)[X]$, it follows from (2.), with $s$ replaced by $s-1$, that $\Phi^{\sp}_{J,(A/PA)^{s-1}}(X) \bmod P \in (C\otimes_A \BF_P)[X]$. Indeed, if this where not the case, consider its highest coefficient not in $C\otimes_A \BF_P$ and remember that all our polynomials are monic.

Lastly, it follows from (1.) 
that now also $\Phi^{\ord}_{J,(A/PA)^s}(X) \bmod P \in (C\otimes_A \BF_P)[X]$.
\end{Proof}

\begin{Ques}
Suppose that $J$ distiguishes the reduced Drinfeld modules $\bar\phi^{(f)}$ for special isogenies $f\in I_P$.
Is $\Phi_{J,(A/PA)^s}^\sp(X)$ irreducible modulo $P$?
\end{Ques}

\section{Examples}

\paragraph{Example 1. $\mathbf s=1$:} In this case 
\[
\Phi_{J,(A/PA)}^\sp(X) = X-J(\phi^{(f_0)}) = X-J^{|P|},
\] 
by Corollary \ref{Frobenius}, so
\[
\Phi_{J,(A/PA)}(X) \equiv (X-J^{|P|})\cdot \Phi_{J,(A/PA)^2}^\sp(X^{|P|}) \bmod P.
\]

\paragraph{Example 2. $\mathbf s=r-1$:} We have
\[
\Phi_{J,(A/PA)^r}^\sp(X) = X-J(\phi^{(\phi_N)}) = X-J,
\]
so
\[
\Phi_{J,(A/PA)^{r-1}}(X) \equiv  \Phi_{J, (A/PA)^{r-1}}^\sp(X) \cdot \left(X^{|P|}-J\right)^{|P|^{r-2}} \bmod P.
\]

\paragraph{Example 3. $\mathbf r=2$:} This is a combination of examples 1 and 2 above, and gives the classical result (see \cite{Bae}):
\[
\Phi_{J,(A/PA)}(X) \equiv (X-J^{|P|})\cdot (X^{|P|}-J) \bmod P.
\]

\paragraph{Example 4.} Now suppose that $r=3$, $P=T$ and $A=\BF_2[T]$, see \cite{BR09}. In this case,
\[
C = \frac{A[J_{07},J_{12}, J_{41}, J_{70}]}{\langle J_{07}J_{41} - J_{12}^4,\quad J_{12}J_{70} - J_{41}^2\rangle},
\]
where $J_{ij} = g_1^i g_2^j$.
In \cite{BR09} we computed $\Phi_{J,(A/TA)^2}(X)$ for $J\in\{J_{07},J_{12},J_{41},J_{70}\}$ (they are denoted $P_{J,T}(X)$ in that paper).
Reducing these modulo $T$, one obtains, for example,
\[
\Phi_{J_{12},(A/TA)^2}(X) \equiv \Phi^{\sp}_{J_{12},(A/TA)^2}(X)\cdot \left(X^2+J_{12}\right)^2 \bmod T
\]
where
\begin{eqnarray*}
  \Phi^{\sp}_{J_{12},(A/TA)^2}(X) & \equiv &  X^{3} + (J_{07} J_{12} + J_{12}^{3} + J_{70}) X^{2}
   + \; (J_{07} J_{41} + J_{12} J_{41} J_{70} + J_{41} + J_{70}^{2}) X  \\
  & & + \;  (J_{07} J_{12} J_{41} + J_{12}^{2} J_{70} + J_{12} J_{70}^{2} + J_{70}^{3} + J_{70})  \bmod T.
\end{eqnarray*}

A similar computation, carried out with the help of Herinniaina Razafinjatovo, confirms that

\[
\Phi_{J_{12},(A/TA)}(X) \equiv (X+J_{12}^2) \cdot \Phi^{\sp}_{J_{12},(A/TA)^2}(X^2) \bmod T,
\]
with $\Phi^{\sep}_{J_{12},(A/TA)^2}(X)$ as above.

\begin{center}
\rule{8cm}{0.01cm}
\end{center}

\begin{minipage}[t]{8cm}{\small
Department of Mathematical Sciences \\
University of Stellenbosch \\
Stellenbosch, 7600 \\
South Africa \\
fbreuer@sun.ac.za}
\end{minipage}
%
\begin{minipage}[t]{8cm}{\small
Institut f\"ur Mathematik  \\
Universit\"at Kassel,\\
Kassel, 34132 \\
Germany \\
rueck@mathematik.uni-kassel.de}
\end{minipage}

\end{document}